\documentclass[10pt]{article}
\usepackage{latexsym,amsmath,amssymb,amsfonts,verbatim}
\usepackage{amsmath,amscd,chngcntr,hyperref}

\author{%
Junwu Tu \thanks{Mathematics Department,
University of Oregon, OR 97403, USA, {\em e-mail: }{\tt junwut@uoregon.edu}}}
\title{Matrix factorizations via Koszul duality}
\date{}
\DeclareFontFamily{U}{rsf}{}
\DeclareFontShape{U}{rsf}{m}{n}{
  <5> <6> rsfs5 <7> <8> <9> rsfs7 <10-> rsfs10}{}
\DeclareMathAlphabet{\mathscr}{U}{rsf}{m}{n}

\DeclareMathAlphabet{\mathgth}{U}{euf}{m}{n}

\DeclareFontFamily{U}{cyr}{}
\DeclareFontShape{U}{cyr}{m}{n}{
  <5> wncyr5 <6> wncyr6 <7> wncyr7 <8> wncyr8 <9> wncyr9 <10-> wncyr10}{}
\DeclareMathAlphabet{\mathcyr}{U}{cyr}{m}{n}

\input cyracc.def

\makeatletter
\def\operator@font{\sf}
\makeatother

\newcommand{\cA}{{\mathscr A}}

\newcommand{\cC}{{\mathscr C}}
\newcommand{\cE}{{\mathscr E}}

\newcommand{\cO}{{\mathscr O}}

\newcommand{\cD}{{\mathscr D}}

\newcommand{\MF}{{\mathsf{MF}}}

\newcommand{\D}{{\mathsf D}_{\mathsf{coh}}^b}
\newcommand{\PN}{{\mathbf{P}}}

\newcommand{\BM}{{\mathsf{BM}}}

\newcommand{\chk}{{\scriptscriptstyle\vee}}

\newcommand{\stab}{{\mathsf{stab}}}

\newcommand{\oppo}{{\mathsf{op}}}

\DeclareMathOperator{\Spec}{Spec}

\DeclareMathOperator{\End}{End}
\DeclareMathOperator{\Hom}{Hom}

\DeclareMathOperator{\Jac}{Jac}
\DeclareMathOperator{\Proj}{Proj}

\DeclareMathOperator{\Tw}{Tw}

\DeclareMathOperator{\id}{id}

\newcommand{\bS}{\mathbf{S}}
\newcommand{\T}{\mathbb{T}}
\newcommand{\B}{\mathbb{B}}

\DeclareMathOperator{\sym}{sym}

\newcommand{\ra}{\rightarrow}

\newcommand{\C}{\mathbb{C}}

\newcommand{\Z}{\mathbb{Z}}

\newcommand{\remark}{\noindent\textbf{Remark: }}
\newcommand{\proof}{\noindent\textbf{Proof. }}

\newcommand{\gr}{{\mathsf{gr}}}
\newcommand{\cone}{{\mathsf{cone}}}

 \newtheorem{theorem}{Theorem}[section]
 
 \newtheorem{lemma}[theorem]{Lemma}
 \newtheorem{corollary}[theorem]{Corollary}
 \newtheorem{proposition}[theorem]{Proposition}

 \newtheorem{definition-theorem}[theorem]{Definition-Theorem}
 \newtheorem{example}[theorem]{Example}
 
\renewcommand{\phi}{\varphi}

\counterwithout{paragraph}{subsection}
\counterwithin{paragraph}{section}
\numberwithin{equation}{section}

\begin{document}
\maketitle
\begin{abstract}

In this paper we prove a version of curved Koszul duality for $\Z/2\Z$-graded curved coalgebras and their coBar differential graded algebras. A curved version of the homological perturbation lemma is also obtained as a useful technical tool for studying curved (co)algebras and precomplexes.

The results of Koszul duality can be applied to study the category of matrix factorizations $\MF(R,W)$. We show how Dyckerhoff's generating results fit into the framework of curved Koszul duality theory. This enables us to clarify the relationship between the Borel-Moore Hochschild homology of curved (co)algebras and the ordinary Hochschild homology of the category $\MF(R,W)$. Similar results are also obtained in the orbifold case and in the graded case.

\end{abstract}
\section{Introduction}
\label{sec:intro}

\paragraph{Backgrounds and motivations.} Matrix factorizations of an element $W$ in a commutative ring $R=\C[[x_1,\cdots,x_n]]$ were first introduced by Eisenbud~\cite{Eisen} in the study of singularity theory. Recently this theory has received renewed interests largely due to its appearance in Kontsevich's homological mirror symmetry conjecture. Indeed the differential graded category $\MF(R,W)$ of matrix factorizations  is conjecturally mirror to the Fukaya category of a Fano symplectic manifold $M$.

The following fundamental results concerning the structure of the dg category $\MF(R,W)$ were obtained by Dyckerhoff~\cite{Dyck} under the assumption that $W$ has isolated singularities:
\begin{description}
\item[---] The homomtopy category $[\MF(R,W)]$ is classically generated by a single object $k^\stab$;
\item[---] The dg algebra $A:=\End_{\MF(R,W)}(k^\stab)$ realizes $\MF(R,W)$ as the dg category of perfect dg modules over $A$;
\item[---] We have $HH_*(\MF(R,W))\cong \Jac(W)[\dim R]$.
\end{description}
Dyckerhoff's computation of the Hochschild homology was indirect, which uses the fact that $\MF(R,W)$ is a compact smooth Calabi-Yau category to reduce the computation to that of Hochschild cohomology. The computation of Hochschild cohomology in turn relies on Toen's interpretation of it as natural transformations from the identity functor to itself.

There is a different approach to understand the category $\MF(R,W)$ initiated in~\cite{Segal} and~\cite{CT}. More precisely the data $(R,W)$ naturally give rise to a curved algebra which we will denote by $R_W$. The category $\MF(R,W)$ can be interpreted as the category of perfect modules over this curved algebra $R_W$. Through this perspective C\u ald\u araru and the author~\cite{CT} introduced the notion of Borel-Moore Hochschild homology of a curved algebra, and proved that
\[ HH_*^\BM(R_W)\cong \Jac(W)[\dim R]\footnote{This isomorphism and the isomorphism below hold for arbitrary smooth commutative ring $R$ and $W$ with isolated singularities.}.\]
Moreover in~\cite{Segal} it was shown that
\[HH_*^\BM(\MF(R,W))\cong HH_*^\BM(R_W).\]
The main advantage of this approach is we have an explicit complex: the Borel-Moore Hochschild chain complex of $R_W$. Thus it would be desirable to relate $HH_*(\MF(R,W))$ with $HH_*^\BM(R_W)$ or $HH_*^\BM(\MF(R,W))$, which would also yield an easier way of computing $HH_*(\MF(R,W))$. To clarify this relationship is the main motivation for the current paper. In the following we explain the main ideas, and give a section-wise summary of our main results.

\paragraph{Curved Koszul duality over $\Z/2\Z$.} The main idea of relating the two types of Hochschild homology is to use Koszul duality theory. In fact Dyckerhoff's results mentioned above already suggest such a link.

For applications to matrix factorizations we need a version of curved Koszul duality theory in the $\Z/2\Z$-graded situation. This theory was developed by Positselski in~\cite{Pos} in great generality where various types of non-standard derived categories were introduced in order to obtain the desired results. In Section~\ref{sec:koszul} we give a more direct proof of the version of Koszul duality which is enough for the applications we have in mind. The proofs we give rely on a curved version of homological perturbation lemma (see Appendix~\ref{sec:chpl}) which is of independent interest. More precisely the main result we obtain in Koszul duality theory is the following theorem.

\begin{theorem}
\label{thm:kos}
Let $B_M$ be a coaugumented curved coalgebra, and let $\Omega B_M$ be its cobar dg algebra. Then there is a quasi-equivalence
\[ \Tw(B_M) \cong \Tw(\Omega B_M)\]
of dg categories of twisted complexes. If furthermore the coalgebra $B$ is conilpotent, then the dg algebra $\Omega B_M$ itself is a compact generator for the homotopy category $[\Tw(\Omega B_M)]$.
\end{theorem}

\medskip
\remark Twisted complexes used in this theorem are not the standard ones in the sense that we allow possibly infinite rank ones, and moreover, we do not assume the upper-triangular condition. Due to these two non-standard conventions, the second part of the above theorem is not at all obvious. We also remark that these modifications are necessary for the purpose of doing Koszul duality.

\paragraph{Applications to $\MF(R,W)$.} In Section~\ref{sec:generator} and Section~\ref{sec:hoch} we apply Koszul duality theory to study $\MF(R,W)$. For this observe that the commutative ring $R$ is the dual algebra of the symmetric coalgebra $C$ generated by variables $y_1:=x_1^\chk,\cdots, y_n:=x_n^\chk$, and the element $W\in R$ is the dual of a linear map $M:C\ra \C$. Moreover matrix factorizations of $(R,W)$ can be identified with twisted complexes over the curved coalgebra $C_M$ which are of finite rank. This simple observation allows us to apply Theorem~\ref{thm:kos} in the situation $C_M$ and $\Omega C_M$ to obtain understanding of $\MF(R,W)$. We summarize our main results in the following theorem.

\begin{theorem}
\label{intro:mf}
Assume that $W$ has isolated singularities. Then we have
\begin{description}
\item[---] Dyckerhoff's generator $k^\stab$ arises from Koszul duality\footnote{For a more precise statement we refer to Section~\ref{sec:generator}.};
\item[---] The dg algebra $A:=\End_{\MF(R,W)}(k^\stab)$ is quasi-isomorphic to the cobar dg algebra $\Omega C_M$;
\item[---] There is a canonical isomorphism $HH_*(\MF(R,W))\cong [HH_*^\BM(R_W)]^\chk$.
\end{description}
\end{theorem}

\medskip
\remark It is an interesting puzzle to understand the appearance of dualization in the above isomorphism between the $HH_*^\BM(R_W)$ and $HH_*(\MF(R,W))$. This might be explained by a relationship between Koszul duality and a natural pairing (generalized Mukai-pairing) on the Hochschild homology.

\paragraph{Applications to $\MF_G(R,W)$.} In Section~\ref{sec:orbifold} we generalize the above results and Dyckerhoff's generation result to the orbifold case. The main results are summarized in the following theorem.

\begin{theorem}
\label{intro:orbi}
Assume that $W$ has isolated singularities, and $G$ a finite abelian group acting on $R$ which fixes $W$. Then we have
\begin{description}
\item[---] The homotopy category $[\MF_G(R,W)]$ of equivariant matrix factorizations is classically generated by
\[\left \{ k^\stab \otimes \C_\chi \mbox{ $|$ $\chi$ is a character for the group } G \right \}\]
where $\C_\chi$ denotes the one dimensional representation associated to the character $\chi$. \item[---] The smash product dg algebra $\Omega(C_M)\sharp G$ realizes $\MF_G(R,W)$ as the dg category of perfect dg modules over $\Omega(C_M)\sharp G$. 
\item[---] For the Hochschild homology we have 
\[HH_*(\MF_G(R,W)) \cong [HH_*^\BM(R_W\sharp G)]^\chk.\]
\end{description}
\end{theorem}

\medskip
\remark In~\cite{CT} the vector space $HH_*^\BM(R_W\sharp G)$ was explicitly computed as
\[ HH_*^\BM(R_W\sharp G) = (\oplus_{g\in G} HH_*^\BM (R_W|_g) )^G \]
where $R_W|_g$ denotes the curved algebra associated to the LG model on the $g$-fixed points of $\Spec(R)$.

\paragraph{Applications to $\MF^\gr(S,W)$.} In Section~\ref{sec:graded} we study graded matrix factorizations where similar results are obtained. In the graded case we consider $S:=\C[x_1,\cdots,x_n]$ the polynomial ring in $n$ variables endowed with its standard grading. Let $W\in S$ be a homogeneous polynomial of degree $d$. Denote $G:=\Z/d\Z$ which acts on $S$ in the obvious way. Since $W$ is of degree $d$ this action preserves $W$. In this situation we can consider the dg category $\MF^\gr(S,W)$ of graded matrix factorizations (see Section~\ref{sec:graded} for its definition).

\begin{theorem}
\label{intro:gradedmf}
Assume that $W$ has isolated singularities. Then we have
\begin{description}
\item[---] The homotopy category $[\MF^\gr(S,W)]$ is classically generated by
\[k^\stab(d-1),k^\stab(d-2),\cdots,k^\stab\]
where the shifts in the parentheses are polynomial degree shifts of graded $S$-modules. 
\item[---] There is a $\Z$-graded~\footnote{This $\Z$-graded is not the standard polynomial grading, see Section~\ref{sec:graded} for its precise definition.} smash product algebra $\Omega(C_M)\sharp G$ realizing $\MF^\gr(S,W)$ as the dg category of perfect dg modules. 
\item[---] For the Hochschild homology we have
\[HH_*(\MF^\gr(S,W))\cong 
[HH_*^\BM (S_W\sharp G)]^\chk\] 
where the operation $\chk$ denotes graded dualization.
\end{description}
\end{theorem}

\paragraph{Acknowledgments.} I would like to thank my advisor Andrei C\u ald\u araru for his continuous support and valuable discussions as well as for reading the first manuscript of this work. I thank Tony Pantev for his encouragement and Bernhard Keller for answering several questions. I also thank Damien Calaque for explaining Koszul duality and Paul Seidel for sharing his unpublished notes. Furthermore I am thankful to Tobias Dyckerhoff and Daniel Pomerleano for pointing out a mistake in an earlier version of the paper and also for making interesting comments. Last but not least I thank Leonid Positselski for answering numerous questions I had in understanding his curved Koszul duality paper. The work is part of the author's thesis which was done in the University of Wisconsin, Madison.

\section{Koszul duality for the (co)bar constructions}
\label{sec:koszul}

In this section we recall the coBar construction and prove a version of Koszul duality between $\Z/2\Z$-graded curved coalgebras and their cobar algebras. Then we prove some useful properties concerning the derived categories of cobar algebras. The results proved in this section are essentially due to Positselski~\cite{Pos}, although we give more direct proofs here.

Throughout this section we will work over a base field $k$. Linear algebra operations such as tensor product or homomorphism between vector spaces are all taken over $k$ unless otherwise stated.
\paragraph{Curved algebras.} A curved algebra structure on a super vector space $A$ is an associative algebra structure on $A$ together with an odd linear map $d:A\rightarrow A$ and an even central element $W\in A$.

\begin{example}
\label{ex:ca}
An example of a curved algebra that will be of primary interest in this paper. Let $V$ be a finite dimensional vector space over a field $k$. Consider the commutative algebra $R:=\widehat{\sym(V^\chk)}$ together with a choice of an element $W$ in it. Since $R$ is commutative any element in it is automatically central.
\end{example}

\paragraph{Twisted complexes over $A_W$: matrix factorizations.} We can define the category $\Tw(A_W)$ of twisted complexes over a curved algebra $A_W$. The objects of this category are pairs $(E,Q)$ where $E$ is a $\Z/2\Z$-graded free $A$-module and $Q$ is an odd $A$-linear map such that $Q^2=W\id$. The morphism space between two objects $(E,Q)$ and $(F,P)$ consists of all $A$-linear maps from $E$ to $F$. As such, the $\Hom$ space inherits a differential defined by $D(\phi)=P\circ \phi -(-1)^{|\phi|} \phi \circ Q$. One easily checks that $D$ squares to zero as $W\id$ is in the center of matrix algebras.

This differential makes the category $\Tw(A_W)$ into a differential graded category. Note that here we allow possibly infinite rank modules in the construction of $\Tw(A_W)$. We denote by $\Tw^b(A_W)$ the full subcategory of $\Tw(A_W)$ consisting of twisted complexes that are of finite rank. The category $\Tw^b(A_W)$ ($\Tw(A_W)$ respectively) is sometimes also referred to the category of matrix factorizations $\MF(A,W)$ ($\MF^\infty(A,W)$ respectively).

As the category $\Tw(A_W)$ has a dg structure we can define the notion of homotopy between morphisms and objects. More precisely, we say two morphisms $f$ and $g$ are homotopic if $f-g$ is exact. We say two objects $E$ and $F$ are homotopic if there are morphisms $f:E\rightarrow F$ and $g:F\rightarrow E$ such that $f\circ g$ is homotopic to $\id_F$ and $g\circ f$ is homotopic to $\id_E$.

\paragraph{Curved coalgebras.} Dualizing the definition for curved algebras we arrive at the definition for curved coalgebras. Namely curved coalgebra structure on a vector space $B$ is a $\Z/2\Z$-graded coassociative coalgebra structure on $B$ together with an even map $M:B\rightarrow k$ such that the composition
\[ B\stackrel{\Delta}{\ra} B\otimes B \stackrel{M\otimes \id-\id\otimes M}{\longrightarrow} B\]
is zero. Using Sweedler's notation for the coproduct $\Delta$, the above is equivalent to
\[M(x^{(1)})x^{(2)}-x^{(1)}M(x^{(2)})=0\]
for all $x\in B$. As before we denote a curved coalgebra by $B_M$. 

\begin{example}
\label{ex:cca}
As a dual example of \ref{ex:ca} we consider $C:=\sym(V)$ be the vector space of symmetric tensors on $V$. Again we consider $C$ as a super vector space concentrated in the even part. There is a natural coalgebra structure on $C=\sym(V)$ whose dual algebra~\footnote{Note that the dual of a coalgebra is always an algebra, but not vice versa due to infinite dimensionality.} is the commutative algebra $R$ in~\ref{ex:ca}. The curvature term is any linear map $M:C\rightarrow k$.
\end{example}

\paragraph{Basics of comodules.} We recall some useful properties of cofree comodules. First of all for purposes of this paper we consider cofree comodules of the form $B\otimes V$ for some $k$-vector space $V$ (possibly infinite dimensional) in the twist construction above. Moreover in the abelian category $\cA$ of $B$-comodules, cofree comodules are injective objects and hence is closed under direct product in $\cA$. For example we have $\prod (B\otimes V_i) \cong B \otimes (\prod V_i)$ where the product on the left hand side is \emph{taken in} $\cA$. A special property for $\cA$ is that the class of injective objects is also closed under direct sum in $\cA$. Explicitly we have $\coprod (B\otimes V_i) \cong B \otimes (\coprod V_i)$. 

\paragraph{Twistec complexes over $B_M$: matrix cofactorizations.} Given a curved coalgebra $B_M$ we can construct a category $\Tw(B_M)$ of twisted complexes. The objects are pairs $(E,Q)$ with $E$ a cofree $B$-comodule and $Q$ an odd comodule map on $E$ such that the dual of the matrix factorization identity holds,
 \[Q^2(x)=M(b^{(1)})x^{(2)}.\]
Here we write the coaction map to be $\rho(x)= b^{(1)}\otimes x^{(2)}$ for $x\in E$, $b^{(1)}\in B$, and $x^{(2)}\in E$. The $\Hom$ spaces and differentials on $\Hom$ spaces are defined in a similar way as for matrix factorizations. Objects in $\Tw(B_M)$ will be called matrix cofactorizations. There is also a differential graded structure on $\Tw(B_M)$. The full subcategory of $\Tw(B_M)$ consisting of matrix cofactorizations that are of finite rank over $B$ will be denoted by $\Tw^b(B_M)$.

There is a simple relation between the two dg categories $\Tw^b(R_W)$ and $\Tw^b(C_M)$, made precise in the following lemma.

\begin{lemma}
\label{lem:dual}
Let $B_M$ be a curved coalgebra, and let $A_W$ be its dual curved algebra. Define a functor $D:\Tw^b(B_M)^\oppo \rightarrow \Tw^b(A_W)$ by formula
\[ (E,Q) \stackrel{D}{\mapsto} (E^\chk, Q^\chk)\]
on objects and for any morphism $f\in\Hom_{\Tw(C_M)}((E,Q),(F,P))$
\[ D(f):= f^\chk:(F^\chk,P^\chk)\rightarrow (E^\chk,Q^\chk). \]
Then $D$ is an equivalence between $\Tw^b(B_M)^\oppo$ and $\Tw^b(A_W)$.
\end{lemma}

\proof Observe that a map $h:B\rightarrow B$ of $B$-comodules is uniquely determined by its composition with the counit map. Conversely any $k$-linear map $\alpha:B \rightarrow k$ defines a map of $B$-comodules by
\[ B\rightarrow B\otimes B \stackrel{\alpha\otimes \id}{\rightarrow} k\otimes B=B.\]
This defines an isomorphism between $\Hom_B(B,B)$ and $\Hom_k(B,k)=B$. More generally for two cofree $C$-comodules $E_1=B\otimes V_1$ and $E_2=B\otimes V_2$ with $V_1$ and $V_2$ finite dimensional vector spaces over $k$ we have
\begin{align*}
 \Hom_B(E_1,E_2)&=\Hom_B(B\otimes V_1,B\otimes V_2) \\
&\cong \Hom_B(B,B)\otimes \Hom_k(V_1,V_2) \cong A \otimes \Hom(V_1,V_2).\end{align*}
For the $\Hom$ space between $DE_2$ and $DE_1$, we have
\begin{align*}
\Hom_A( (B\otimes V_2)^\chk, (B\otimes V_1)^\chk)&= \Hom_A(A\otimes V_2^\chk, A\otimes V_1^\chk)\\
& = A\otimes \Hom(V_2^\chk, V_1^\chk)= A\otimes \Hom(V_1,V_2)\end{align*}
where the first and the last equality follow from $V_1$ and $V_2$ being finite dimensional. Thus we have verified that the functor $D$ is an equivalence. A direct computation shows that it also preserves the differential and hence the lemma is proved.

\paragraph{The coBar construction.} Let $B_M$ be a curved coalgebra. A $k$-linear map $\eta:k\rightarrow B$ is called a coaugumentation of $B_M$ if 
\begin{itemize}
\item[(1)] $\eta$ splits the counit map;
\item[(2)] $\eta$ is a map of coalgebras;
\item[(3)] $M\circ \eta=0$.
\end{itemize}
Denote by $B^+$ the cokernel of $\eta$ which can be identified with the kernel of the counit through the splitting in $(1)$. Indeed the above conditions imply there is a direct sum decomposition $B\cong B^+\oplus k$ of $B$ as coalgebras, moreover $M$ vanishes on the component $k$.

Given a coaugumented curved coalgebra $B_M$, we can construct a differential graded algebra  $\Omega B_M$, known as its coBar construction. Explicitly as an associative algebra $\Omega B_M$ is the free tensor algebra generated by $B^+[-1]$ which is simply
\[\T B^+[-1]= \oplus_{k=0}^\infty (B^+[-1])^{\otimes k}.\]
The differential $d$ is a derivation on $\Omega B_M$ determined by the following two components
\begin{align*}
B^+\hookrightarrow B\rightarrow &B\otimes B \rightarrow B^+\otimes B^+;\\
B^+\hookrightarrow B & \stackrel{M}{\rightarrow} k.
\end{align*}

\begin{example}
\label{ex:cobar}
Let us work out the coBar construction of the curved coalgebra $C_M$ defined in Example~\ref{ex:cca}. There is a natural coaugmentation on $C$: the inclusion of scalars. In order for it be compatible with curvature, we assume that $M$ vanishes on scalar part of $C$.This coaugmentation induces in particular a direct sum decomposition $C\cong C^+\oplus k$. The coBar algebra $\Omega C_M$ is the free tensor algebra generated by $\sym(V)^+[-1]$ with differential given by the sum of two components which we denote by $d^+$ and $d^-$. These maps act on monomials $f_1|\cdots|f_k$ by
\begin{align*}
d^+(f_1|f_2|\cdots|f_k)&=\sum_{i=1}^k (-1)^{i-1} f_1|\cdots|\Delta(f_i)|\cdots|f_k,\\
d^-(f_1|f_2|\cdots|f_k)&=\sum_{i=1}^k (-1)^{i-1} M(f_i)f_1|\cdots|\widehat{f_i}|\cdots|f_k
\end{align*}
where $\Delta$ is the coproduct on $C^+$ induced from that of $C$.
\end{example}

\paragraph{Twisting cochains.} For a curved coalgebra $B_M$ and a unital differential graded algebra $A$, one can construct a curved differential graded algebra structure on the space of $k$-linear maps $\Hom(B,A)$. It is defined by the following formulas:
\begin{itemize}
\item Curvature: $W(B,A)$: $B\stackrel{M}{\rightarrow}k\stackrel{\mbox{unit}}{\rightarrow} A$;
\item Differential: $(d\phi)(x)= d(\phi(x))$;
\item Product: $(\phi*\psi)(x)= (-1)^{|x^{(1)}||\psi|}\phi(x^{(1)})\psi(x^{(2)})$.
\end{itemize}
A twisting cochain from $B$ to $A$ is an odd element $\tau\in \Hom(B,A)$ such that
\[ \tau*\tau +d\tau + W(B,A)=0. \]
There is a natural twisting cochain from $B_M$ to its coBar algebra $\Omega B_M$ defined by
\[ \tau_{B_M}: B \rightarrow B^+\hookrightarrow \Omega B_M.\]
\paragraph{Correspondence of twisted complexes.} We can use the twisting cochain $\tau_{B_M}$ to define a correspondence between categories of twisted complexes. We work out explicitly this correspondence for a coaugumented curved coalgebra $B_M$ its coBar algebra $\Omega B_M$. The goal is to construct dg functors 
\begin{align*}
\Phi:& \Tw(B_M) \rightarrow \Tw(\Omega B_M)\\
 \Psi:&\Tw(\Omega B_M)\rightarrow \Tw(B_M).\end{align*}
We begin with the construction of $\Phi$. Let $(E,Q)$ be a matrix cofactorization over $B_M$. We need produce a twisted complex $\Phi(E)$ over $\Omega B_M$. As a vector space over $k$ it is simply $\Omega B_M\otimes E$. The left $\Omega B_M$-module structure is induced from that of $\Omega B_M$. The differential on $\Omega B_M\otimes E$ is defined using the natural twisting cochain $\tau_{B_M}$:
\[ d(x\otimes e)=dx\otimes e +(-1)^{|x|}x\otimes Qe + (-1)^{|x|+1}x \tau(y^{(1)})\otimes e^{(2)}\]
where we have denoted the coaction map $\rho:E\rightarrow B\otimes E$ by $\rho(e)=y^{(1)}\otimes e^{(2)}$ for $y^{(1)}\in B$. One checks that $d^2=0$ and that it is compatible with the left module structure on $\Phi(E)$. We write $\Phi(E)=\Omega\otimes^\tau E$ where the superscript $\tau$ is to indicate that we are using the twisting cochain $\tau$ to define the differential on $\Phi(E)$. Note that $\Phi(E)$ is of infinite rank whenever $B$ is of infinite dimension over $k$. For this reason we need to consider $\Tw(\Omega B_M)$ instead of $\Tw^b(\Omega B_M)$. For a morphism $f:(E,Q)\rightarrow (F,P)$ in $\Tw(B_M)$, define $\Phi(f)=\id\otimes f$ from $\Phi(E)$ to $\Phi(F)$. One can check that $\Phi$ is a differential graded functor between differential graded categories.

In the reverse direction, if $(F,d)$ is a twisted complex over $\Omega B_M$, we need to define a matrix cofactorization $\Psi(F)$ over $B_M$. As a vector space this is $B\otimes F$. The left $B$-comodule structure is induced from that of $B$ and the matrix cofactorization map is defined by
\[ Q(x\otimes f)= dx\otimes f +(-1)^{|x|}x\otimes df +x^{(1)}\otimes \tau(x^{(2)}) f\]
where $\tau(x^{(2)})f$ is the action of $\Omega B_M$ on $F$. One checks that $Q$ satisfies the matrix cofactorization identity and hence defines a twisted complex (again of infinite rank) over $B_M$. Similarly the above construction extends to the morphism space and hence defines a dg functor $\Psi$ in the reverse direction.

\paragraph{Curved coBar duality over $\Z/2\Z$.} As both the categories $\Tw(B_M)$ and $\Tw(\Omega B_M)$ are dg categories, one can speak of the notion of homotopy between dg functors. Namely we say a functor $F$ is homotopic to another functor $G$ if they are homotopic when applied to any object. Hence we can also define the notion of homotopy between categories by requiring that there are functors in both ways such that their compositions are homotopic to the identity functors. The following theorem is the Koszul duality property for the cobar construction. Essentially it is duality between the categories of twisted complexes.

\begin{theorem}
\label{thm:koszul}
The functors $\Phi$ and $\Psi$ are homotopy inverse of each other. Hence the two categories $\Tw(B_M)$ and $\Tw(\Omega B_M)$ are homotopic.
\end{theorem}

\proof We start by showing that the composition $\Psi\circ \Phi$ is homotopic to the identity functor on $\Tw(B_M)$. For any object $(E,Q)\in \Tw(B_M)$, consider the morphism $\eta_E$ between $E$ and $\Psi\Phi(E)=B\otimes^\tau \Omega B_M \otimes^\tau E$ defined by 
\[ E\rightarrow B\otimes E \hookrightarrow B\otimes^\tau \Omega B_M \otimes^\tau E\]
where the first map is the coaction map and the second map is simply putting unit of $\Omega B_M$ on the middle position of $B\otimes^\tau\Omega B_M\otimes^\tau E$. A direct computation shows that $\eta_E$ is a map of twisted complexes over $B_M$. 
\begin{lemma}
\label{lem:cone}
Let $f:(E,Q)\rightarrow (F,P)$ be a closed morphism in $\Tw(B_M)$. Define the cone of $f$ to be the matrix cofactorization $(E[1]\oplus F,T)$ with $T$ given by the matrix
\[ T=\left[ \begin{array}{cc}
Q & 0  \\
f & P  \end{array} \right].\]
Then $f$ is a homotopy equivalence if and only if $\cone(f)$ is contractible. 
\end{lemma}

\proof If $\cone(f)$ is contractible, there exists a morphism $H:\cone(f) \rightarrow \cone(f)$ such that
\[ \id =[T,H] .\]
Writing $H$ as a matrix
\[ \left[ \begin{array}{cc}
a & b  \\
c & d  \end{array} \right],\]
after a matrix multiplication, we find that the map $b$ defines a homotopy inverse of $f$. A similar consideration works for the reversed direction. The lemma is proved.

\medskip
\noindent By the above lemma, in order to show that $E$ and $\Psi\Phi(E)$ are homotopic it is enough to prove that $\cone(\eta_E)$ is contractible. The cone $\cone(\eta_E)$ is explicitly given by $E[1]\oplus (B\otimes ^\tau\Omega B_M\otimes ^\tau E)$ on which acts an operator $D$ which satisfies the matrix cofactorization identity. We next write down the map $D$ on $E[1]$ and $(B\otimes ^\tau\Omega B_M\otimes ^\tau E)$.
\begin{comment}
This space has a grading by the number of $B$-tensors. We write $b_0[b_1|\cdots|b_l]\otimes e$ for an element containing $l+1$ $B$-tensors and $e$ for an element with no $B$-tensors .\end{comment}
On elements of the form $b_0[b_1|\cdots|b_l]\otimes e\in (B\otimes ^\tau\Omega B_M \otimes ^\tau E)$ the predifferential $D$ acts by
\begin{align*}
D&:=d_\Delta + Q + d_M;\\
d_\Delta (b_0[b_1|\cdots|b_l]\otimes e)&:= \sum_{i=0}^{l} (-1)^i b_0[b_1|\cdots|\Delta(b_i)|\cdots|b_l ]\otimes e +\\
&+b_0[b_1|\cdots|b_l|b^{(1)}]\otimes e^{(2)};\\
Q(b_0[b_1|\cdots|b_l]\otimes e)&:= (-1)^{l+1}b_0[b_1|\cdots|b_l]\otimes Qe;\\
d_M (b_0[b_1|\cdots|b_l]\otimes e)&:= \sum_{i=1}^{l} (-1)^i b_0[b_1|\cdots|M(b_i)|\cdots|b_l]\otimes e.
\end{align*}
On elements in $E[1]$ the predifferential $D$ acts by
\begin{align*}
D&:=d_\Delta + Q + d_M;\\
d_\Delta (e) &:= b^{(1)}\otimes e^{(2)};\\
Q(e)&:= Q(e);\\
d_M(e)&:=0.
\end{align*}
We observe that the differential $d_\Delta$ is simply the cobar resolution of the $B$-comodule $E$
\[ E\rightarrow B\otimes E \rightarrow B\otimes B^+ \otimes E \rightarrow \cdots\]
which is exact. Moreover since the $B$-comodule $E$ is cofree (hence injective) there exists a $B$-linear homotopy $H$ on the cobar resolution above that makes the complex contractible over $B$. Note that the homotopy reduces the number of $B$-tensor components by one. This homotopy operator $H$ defines a homotopy retraction data $(0, 0, H)$ between the zero complex and the cobar resolution (see the Appendix~\ref{sec:chpl} for details on homological perturbation technique). We also want to require $H$ to be special, i.e. $H^2=0$. This can be achieved by making the following transformation
\[ H\mapsto Hd_\Delta H.\]
As the maps $d_\Delta$ are also $B$-linear, the special homotopy retraction is also $B$-linear. To show that $\cone(\eta_E)$ is contractible, we need to show that there exists a $B$-linear homotopy for $D$. For this we consider $D$ as obtained from $d_\Delta$ by a small perturbation $Q+d_M$. Then apply homological perturbation lemma to obtain the homotopy for $D$. As mentioned earlier, the map $D$ is not really a differential as $D^2$ is not zero. Thus the ordinary homological perturbation lemma does not apply to this case. However $D$ satisfies the matrix cofactorization identity by its construction. In this situation a curved version of  homological perturbation lemma can still be applied as is explained in the Appendix~\ref{sec:chpl}. To perform perturbation we need the following lemma.

\begin{lemma}
\label{lem:small}
The curved perturbation $\delta:=Q+d_M$ is small. That is we can define the operator $(\id-\delta\circ H)^{-1}$ on $\cone(\eta_E)$. In fact the operator $\delta\circ H$ is locally nilpotent on $\cone(\eta_E)$.
\end{lemma}

\proof For a $\Z^{\geq 0}$-graded vector space we say an operator on it is locally nilpotent if for any element of bounded degree it is nilpotent. In our case, we consider the space $\cone(\eta_E)$ be graded by the number of $B$-tensors. Observe that the operator $Q$ preserves the number of $B$-tensors while $d_M$ reduces the number of $B$-tensors by one. The homotopy operator also reduces the number of $B$-tensors by one. Hence the composition $\delta\circ H$ strictly reduces the number of $B$-tensors. So it must be locally nilpotent by degree consideration. Since $\delta\circ H$ is a locally nilpotent operator, one can define the operator $(\id-\delta\circ H)^{-1}$ on the direct sum of each graded components which is $\cone(\eta_E)$~\footnote{Note that it is important here that here $\cone(\eta_E)$ is a direct sum rather than a direct product.}. The lemma is proved.

\medskip
\noindent
Applying the curved homological perturbation lemma~\ref{lem:chpl} over the linear category of $B$-linear morphisms, we conclude that there exists a homotopy $H_1$ for the operator $D$.  Hence $\cone(\eta_E)$ is contractible. We have finished half of the proof of the theorem, namely the composition $\Psi\Phi$ is homotopic to the identity functor on $\Tw(B_M)$.

The other half of the proof, the fact that $\Phi\Psi$ is homotopic to identity on $\Tw(\Omega B_M)$ is relatively easier as we are dealing with actual complexes. The author learned this argument from Positselski. It is an expanded version of the proof given in subsection $6.4$ in~\cite{Pos}. For an object $F\in \Tw(\Omega B_M)$ we consider the natural map $\epsilon_F: \Phi\Psi (F):= \Omega B_M \otimes^\tau B_M \otimes^\tau F \rightarrow F$ defined by
\[  a\otimes 1\otimes f \mapsto a f \]
and zero on the other tensors. To show that $\epsilon_F$ is a homotopy equivalence it suffice to show that the cocone $K:=\cone(\epsilon_F)[-1]$ is contractible. This dg module as a vector space is $F[-1]\oplus \Omega B_M\otimes B \otimes F$. Define a finite decreasing filtration on it by
\[ F^0K:=K \supset F^1K:= \Omega(B_M)\otimes B\otimes F \supset F^2K:= \Omega(B_M)\otimes B^+ \otimes F \supset F^3K:=0.\]
One checks that the differential on $K$ does not preserve this filtration but sends $F^iK$ to $F^{i-1}K$. Moreover the induced differential on the associated graded components agrees with the canonical resolution
\[ 0\rightarrow \Omega B_M\otimes B^+\otimes F \rightarrow \Omega B_M\otimes k \otimes F \rightarrow F \rightarrow 0\]
which is exact. Then we can define a dg $\Omega B_M$-submodule of $K$ by 
\[ L:= F^2K + d F^2K\]
where $d$ is the differential on $K$. It follows from the exactness of the above short exact sequence that both $L$ and $K/L$ are contractible. In general this does not imply that $K$ is also contractible. But in our case the dg module $K/L$ is free as $\Omega B_M$-modules, which implies that $K$ admits a direct sum decomposition $L\oplus K/L$ as $\Omega B_M$-modules. Note that this splitting does not necessarily preserve the differential on $K$, nevertheless it realizes $K$ as the cone of a closed map from $L[-1]$ to $K/L$, which implies that $K$ itself is also contractible. The proof of Theorem~\ref{thm:koszul} is now complete.

\paragraph{Homological properties of $\Tw(\Omega B_M)$.} We first introduce some notations. For a dg category $\cD$ we denote by $[\cD]$ its homotopy category. Recall that $[\cD]$ has the same objects as $\cD$, but the morphism spaces between objects are given by the zeroth cohomology of the morphism spaces in $\cD$. Our next goal is to have an understanding of the category $[\Tw(\Omega B_M)]$.

It is a well-known fact that for a dg algebra $A$ the category $[\Tw(A)]$ is a triangulated category. However it \emph{does not} agree with the derived category of $A$ in general. The reason is that the derived category of $A$ is defined by the localization of $[\Tw(A)]$ with respect to the class of acyclic objects (dg modules with zero cohomology) which might not be trivial in $[\Tw(A)]$. Equivalently this is to say that there might exist objects in $\Tw(A)$ that are acyclic while not being contractible. One such example is to take $A=k[x]/x^2$ and $E\in \Tw(A)$ to be
\[ \cdots A \rightarrow A \stackrel{\cdot x}{\rightarrow} A \rightarrow A \cdots\]
where the maps are all given by multiplication by $x$. Then $E$ is acyclic while it is not contractible.

However for a coaugumented conilpotent coalgebra $B$ endowed with a curvature term $M$, we will show that acyclic objects are the same as contractible objects in $\Tw(\Omega B_M)$. Recall that a coaugumented coalgebra $B$ is conilpotent if $B^+$ is the union of the kernels of finite iterated coproducts.

\begin{proposition}
\label{prop:acyclic}
Le $B$ be a coaugumented conilpotent coalgebra and let $F$ be an object in $\Tw(\Omega B_M)$. Then $F$ is acyclic if and only if $F$ is contractible.
\end{proposition}

\proof It suffices to prove that if $F$ is acyclic then it is contractible. As $F$ is an acyclic complex there always exists a contracting homotopy for $F$ over the field $k$. Let $H$ be such a $k$-linear special homotopy of $F$. Consider the Koszul dual $\Psi(F)=B\otimes ^\tau F$. The $B$-linear map $\id\otimes H$ defines a special contracting homotopy for the complex $(B\otimes F, \id \otimes d_F)$. The predifferential $Q$ on $\Psi(F)$ is given by
\[ Q= \id \otimes d_F + d^\tau\]
where the map $d^\tau$ comes from the natural twisting cochain $\tau$ associated with the curved coalgebra $B_M$. We consider $\delta:=d^\tau$ as a curved perturbation of $\id\otimes d_F$ and apply the curved homological perturbation lemma as in the proof of the Theorem~\ref{thm:koszul}. For this we need to prove the curved perturbation $\delta$ is small.
This is an immediate consequence of the conilpotency condition on $C$. In fact the conilpotency condition implies that $\delta \circ (\id\otimes H)$ is a locally nilpotent operator. Thus by the curved homological perturbation lemma~\ref{lem:chpl} the object $\Psi(F)$ is contractible in $\Tw(B_M)$. It follows that the object $\Phi\Psi(F)$ is also contractible. By Theorem~\ref{thm:koszul}, $\Phi\Psi(F)$ is homotopic to $F$ and hence $F$ is also contractible.
Thus the proposition is proved.

\paragraph{Terminologies about generators.} Proposition~\ref{prop:acyclic} immediately implies that the dg algebra $\Omega B_M$ itself is a generator for the triangulated category $[\Tw(\Omega B_M)]$ if $B$ is conilpotent. To make a more precise statement we recall several distinct notions of generators for triangulated categories. We follow the exposition in~\cite{BV}. Let $\cD$ be a triangulated category. A set of objects $\cE:=\left \{ E_i | i\in I \right \}$ is said to classically generate $\cD$ if the smallest triangulated subcategory of $\cD$ containing $\cE$ that is closed under isomorphism and direct summands is equal to $\cD$ itself. We say that $\cD$ is finitely generated if it is classically generated by one object.

The second notion of generation is defined via the orthogonal category of $\cE$. Namely, we say
that $\cE$ weakly generates $\cD$ if the right orthogonal $\cE^\perp$ is trivial. (The right
orthogonal $\cE^\perp$ is by definition the full subcategory of $\cD$ consisting of
objects $A$ such that $\Hom_\cD(E_i[n],A)=0$ for all $i$ and all $n$.) It is clear that classical generators are also weak generators. But the converse is not true in general, often we will drop the adverb "weak" and say that $\cE$ generates $\cD$ if $\cE$ weakly generates it.

If furthermore the category $\cD$ admits arbitrary direct sums one can define the notion of compactness for objects. In such a category an object $E$ in $\cD$ is
said to be compact if the functor $\Hom_\cD(E,-)$ commutes with direct sums. Denote by
$\cD^c$ the full subcategory consisting of compact objects. We say that $\cD$ is compactly
generated if $\cD^c$ generates $\cD$. We need the following result by Ravenel and Neeman~\cite{Neeman}.
\begin{theorem}
\label{thm:rn}
Assume that a triangulated category $\cD$ admitting arbitrary coproduct is compactly generated. Then a set of compact objects classically generates $\cD^c$ if and only if it generates $\cD$.
\end{theorem} 

\begin{corollary}
\label{cor:generator}
Let the notations and assumptions be the same as in Proposition~\ref{prop:acyclic}. Then the dg-module $\Omega(B_M)$ is a compact generator for the category $[\Tw(\Omega B_M)]$. Moreover it classically generates the compact subcategory $[\Tw(\Omega B_M)]^c$.
\end{corollary}

\proof It is clear that the object $\Omega B_M$ is compact. Moreover if $F\in [\Tw(\Omega B_M)]$ is right orthogonal to $\Omega B_M$, it implies that the object $F$ is acyclic. Then it follows from Proposition~\ref{prop:acyclic} that $F$ is in fact contractible hence becomes zero in $[\Tw(\Omega B_M)]$. The last assertion follows from Theorem~\ref{thm:rn}.

\section{Generators for $\MF(R,W)$}
\label{sec:generator}
In this section we work with the curved coalgebra $C_M$ and its dual curved algebra $R_W$ as introduced in Examples~\ref{ex:ca} and~\ref{ex:cca}. As symmetric coalgebras with their canonical coaugmentations are conilpotent coalgebras, all the results in the previous section hold for $C_M$. We prove that the image of the cobar algebra $\Omega C_M$ itself under the Koszul duality functor lies in $\Tw^b(C_M)$. Hence its $k$-linear dual makes sense and defines a matrix factorization in $\Tw^b(R_W)=\MF(R,W)$. Then we identify it with Dyckerhoff's $k^\stab$. Corollary~\ref{cor:generator} and Proposition~\ref{prop:compact_generator} then implies a homological interpretation for $k^\stab$ to classically generate $\MF(R,W)$. This homological interpretation is used in Sections~\ref{sec:orbifold}, ~\ref{sec:graded} to produce classical generators for the derived categories of equivariant or graded matrix factorizations.

\paragraph{Compact generator for $[\Tw C_M]$.} We begin to construct a compact generator for the homotopy category of $\Tw(C_M)$. Note that it is clear that in both the category $\Tw(C_M)$ and $\Tw(\Omega C_M)$ arbitrary coproducts exist and hence one can talk about compactness of objects in these categories. By Theorem~\ref{thm:koszul} the two dg categories $\Tw(C_M)$ and $\Tw(\Omega C_M)$ are homotopic via the homotopy equivalences $\Phi$ and $\Psi$ that preserve coproducts. Hence $\Phi$ and $\Psi$ send compact generators to compact generators. By Corollary~\ref{cor:generator} the object $\Omega C_M$ is a compact generator for the homotopy category $[\Tw(\Omega C_M)]$ as symmetric coalgebras are conilpotent. It follows that the matrix cofactorization $\Psi(\Omega C_M)$ is a compact generator for the homotopy category of $\Tw(C_M)$.

\begin{proposition}
\label{prop:compact_generator}
The homotopy category of $\Tw(C_M)$ is compactly generated by $\Psi(\Omega C_M)$. Moreover $\Psi(\Omega C_M)$ is homotopic to an object in $\Tw^b(C_M)$.
\end{proposition}

\proof Previous discussions have already proved the first assertion. We only need to prove the second assertion. The idea is again to use the curved homological perturbation lemma~\ref{lem:chpl}. By definition the predifferential $Q$ on $\Psi(\Omega C_M):=C\otimes^\tau \Omega C_M$ can be split into three parts defined by
\begin{align*}
d^+(x\otimes y)&:=x\otimes d^+(y) ;\\
d^-(x\otimes y)&:=x\otimes d^-(y);\\
d^\tau(x\otimes y)&:= x^{(1)}\otimes \tau(x^{(2)})y;\\
Q&:=d^++d^-+d^\tau.
\end{align*}
Consider the sum $\delta:=d^-+d^\tau$ as a curved perturbation for the operator $d^+$. We can choose a $k$-linear special homotopy $H$ (always exists over a field) between $(\wedge^*(V),0)$ and $(\Omega C_M, d^+)$ such that $H$ deceases the tensor degree (as $d^+$ increases it). Then extend it to a special homotopy between
\[ (C\otimes \wedge^*(V),0) \cong (C\otimes \Omega C_M, d^+)\]
by putting $\id$ on the $C$ part. To see that the curved perturbation $\delta:=d^-+d^\tau$ is small, note that $d^-$ reduces the number of tensor components, $d^\tau$ reduces the degree of the $C$ part and $H$ reduces the number of tensor components. This allows us to apply the curved homological perturbation lemma~\ref{lem:chpl} to $\Psi(\Omega C_M)$, which implies that $\Psi(\Omega C_M)$ is homotopic to a matrix cofactorization on $C\otimes \wedge^*(V)$. Thus the proposition is proved. This matrix cofactorization obtained via perturbation will still be denoted by $\Psi(\Omega C_M)$.

\paragraph{Relationship with Dyckerhoff's generator $k^\stab$.} In~\cite{Dyck} Dyckerhoff defined a matrix factorization on $R\otimes \wedge^*(V^\chk)$ which he denoted by $k^\stab$. The space $k^\stab$ is a super space with parity determined by the exterior degree. The matrix factorization on $k^\stab$ is defined by choosing a basis $ x_1,\cdots,x_n $ of $V^\chk$, and write $W$ in the form $\sum_{i=1}^n x_iW_i$. Denote the dual basis for $V$ by $ y_1,\cdots,y_n $. Then the matrix map $Q^\chk$ is defined by
\[ Q^\chk(f\otimes \alpha):= x_i f \otimes \lrcorner y_i \alpha + W_i f\otimes x_i\wedge \alpha\]
where $\lrcorner y_i$ denotes the contraction operator and repeated indices are implicitly summed. The goal here is to compare Dyckerhoff's $k^\stab$ with $D\Psi(\Omega C_M)$ produced by Koszul duality (where $D$ is the dualizing functor between cofactorizations and factorizations). The following proposition proves that they are homotopic objects. This can be viewed as a generalization of the classical fact that $D\Psi(\Omega(C))$ is homotopic to the Koszul complex of the residue field $k$.

\begin{proposition}
\label{prop:same}
With the notations introduced above we have a homotopy equivalence
\[ D(\Psi(\Omega C_M)) \cong k^\stab\]
between objects in $\MF(R,W)$.
\end{proposition}

\proof Since the functor $D$ is an equivalence of categories, we denote by $E:=(E,Q)$ the matrix cofactorization whose dual is $k^\stab$. As $\Psi$ is a homotopy inverse to $\Phi$, it is enough to prove that
\[ \Phi \circ \Psi(\Omega C_M) \cong \Phi(E).\]
As shown in the proof of Theorem~\ref{thm:koszul} the counit of the adjunction map
\[ \Phi\circ \Psi(\Omega C_M) \stackrel{\epsilon_{\Omega C_M}}{\rightarrow} \Omega C_M\]
is a homotopy equivalence. Hence it suffices to show that $\Phi(E)$ and $\Omega C_M$ is homotopic. The object $\Phi(E)$ as a vector space is given by $\Omega C_M\otimes C\otimes \wedge^*(V)$. 
Define a linear map $\alpha$ from $\Omega C_M$ to $\Phi(E)$ by
\[ [f_1|\cdots|f_k]\mapsto [f_1|\cdots|f_k]\otimes 1\otimes 1\]
where the middle $1$ is the image of the coaugmentation map of $1\in k$. The last $1$ is the unit in $\wedge^*(V)$. The map $\alpha$ clear respects the left $\Omega(C_H)$-module structure. Moreover it is a map of complexes as $Q$ vanishes on $1\otimes 1$ ($Q^\chk$ increase the polynomial degree on $C$, $Q$ must decrease the degree). We use homological perturbation to show that $\alpha$ is a homotopy equivalence. Again we split the differential $D$ on $\Phi(E)$ into several parts and use homological perturbation lemma. Explicitly for an element $a\otimes f \otimes y \in \Omega C_M\otimes C \otimes \wedge^*(V)$, the map $D$ is the sum of the following four parts:
\begin{align*}
d_\Omega (a \otimes f \otimes y) &:= d_\Omega (a) \otimes f \otimes y;\\
d^\tau (a\otimes f\otimes y) &:= a \tau(f^{(1)}) \otimes f^{(2)} \otimes y;\\
Q^+ (a\otimes f\otimes y) &:= a \otimes \frac{\partial f}{\partial y_i} \otimes y_i \wedge y;\\
Q^-(a\otimes f\otimes y) &:= a \otimes D_i(f)\otimes \lrcorner x_i
\end{align*}
where $y_i$ as before is a basis for the vector space $V$. The map $D_i$ is defined by
\[ C\rightarrow C\otimes C \stackrel{D(W_i\otimes \id)}{\longrightarrow} k\otimes C=C.\]
The map $Q^+$ is simply the Koszul differential on $C\otimes \wedge^*(V)$. We consider the differential $d:=d_\Omega+Q^+$ on the underlying vector space of $\Phi(E)$ and the other part $\delta:=d^\tau +Q^-$ as perturbations of $d$.  One can easily write down a special homotopy $H$ for the Koszul differential $Q^+$ and extend it by $\id$ on $\Omega C_M$ to give a homotopy retraction data between $\Omega C_M$ and $(\Omega C_M\otimes C\otimes \wedge^*(V), d)$. The fact that the perturbation $\delta$ is small follows from the conilpotency property of $C$ and that the curvature $M$ vanishes on scalar and linear terms. Moreover observe that both $H$ and $\delta$ are $\Omega C_M$-linear and
\[ \delta \circ \alpha =0,\]
which implies that the perturbed inclusion is still $\alpha$ and the perturbed differential is still $d_\Omega$ on $\Omega C_M$ by formulas in Appendix~\ref{sec:chpl}. Hence the proposition is proved.

\remark It follows from this Proposition that the endomorphism dg algebra $\End(k^\stab)$ is homotopic to $\Omega C_M$. One can easily prove that the homology of $\Omega C_M$ is $\wedge^*(V)$ assuming that $W$ vanishes on scalars and linear terms. The minimal model $A_\infty$ algebras on $\End(k^\stab)$ has been studied by Dyckerhoff and here we could use $\Omega C_M$ to obtain similar results.

\paragraph{Discussion on generating results.} The two propositions~\ref{prop:compact_generator}, ~\ref{prop:same} above suggest a new homological proof of the fact that $k^\stab$ is a generator for $[\MF(R,W)]$. Indeed it is a direct consequence of these two propositions that $k^\stab$ weakly generates $[\MF(R,W)]$, i.e. its right orthogonal full subcategory is trivial.

However it does not imply that $k^\stab$ classically generates $[\MF(R,W)]$. The problem here is the subcategory $\Tw^b(C_M)$ might not be compact in $\Tw(C_M)$. Indeed we show that this is equivalent to the condition that the object $k^\stab$ classically generates $[\MF(R,W)]$. We need the following theorem (which can be found in~\cite{Neeman}) that characterizes compact objects.
\begin{theorem}
\label{thm:character}
Let $\cD$ be a triangulated category with arbitrary coproduct. Moreover, assume that $\cD$ is compactly generated by a set of compact objects $\cE$. Then an object of $\cD$ is compact if and only if it is a direct summand of an iterated extension of copies of objects of $\cE$ shifted in both directions.
\end{theorem}

\begin{proposition}
\label{prop:equiv}
The full subcategory $[\Tw^b(C_M)]$ of $[\Tw(C_M)]$ is compact if and only if $k^\stab\cong D\Psi (\Omega C_M)$ classically generates $[\MF(R,W)]$.
\end{proposition}

\proof Assume that $[\Tw^b(C_M)]$ is a compact subcategory of $[\Tw(C_M)]$, i.e. every object of $[\Tw^b(C_M)]$ is compact, then it follows from Theorem~\ref{thm:character} that the object $\Psi(\Omega C_M)$ in $[\Tw^b(C_M)]$ obtained by perturbation classically generates $[\Tw^b(C_M)]$ as it is a compact generator for $[\Tw(C_M)]$. Apply the equivalence functor $D$ implies that $D\Psi(\Omega C_M)=k^\stab$ classically generates $[\MF(R,W)]$.

Conversely, if $k^\stab$ classically generates $[\MF(R,W)]$, by Theorem~\ref{thm:character} we conclude that objects in $[\Tw^b(C_M)]$ can be obtained from $\Psi(\Omega C_M)$ by taking direct factors of iterated extensions and shifts, which implies that objects in $[\Tw^b(C_M)]$ are compact in $[\Tw(C_M)]$ as $\Psi(\Omega C_M)$ is a compact generator.

We will now show that the homological smoothness of the dg algebra $\Omega C_M$ implies that the object $k^\stab$ classically generates $[\MF(R,W)]$. Recall that a dg algebra $A$ is called homologically smooth if $A$ considered as an $A\otimes A$-bimodule is a perfect object, i.e. it is a direct factor of finite rank free $A\otimes A$ dg-module.
\begin{proposition}
\label{prop:smooth}
If the dg algebra $\Omega C_M$ is homologically smooth then the full subcategory $[\Tw^b(C_M)]$ of $[\Tw(C_M)]$ is compact.
\end{proposition}

\proof A matrix cofactorization structure on $C\otimes V$ is equivalent to a $\Omega C_M$ dg-module structure on $V$. Hence it suffices to show that any finite dimensional dg $\Omega C_M$-module is compact in $\Tw(\Omega C_M)$. Homological smoothness implies the existence of resolution of diagonal by a perfect complex of $\Omega C_M\otimes \Omega C_M$-bimodules. Via integral transform it produces a resolution for any finite dimensional dg module by a perfect complex of $\Omega C_M$- modules. Thus the proposition is proved.

\section{Hochschild invariants}
\label{sec:hoch}

As another application of Theorem~\ref{thm:koszul} we show that one can calculate the Hochschild homology of $\MF(R_W)$ using the Borel-Moore Hochschild chain complex of the curved algebra $R_W$. The latter was introduced and explicitly computed in~\cite{CT}. We assume that $W$ has isolated singularities throughout this section.

\paragraph{Reducing to Hochschild homology of $\Omega C_M$.} As mentioned in Section~\ref{sec:koszul} the dg category $\Tw^b(R_W)$ is isomorphic as a dg category to $\Tw^b(C_M)^{op}$. Since the Hochschild homologies for opposite dg categories are isomorphic, we have
\[ HH_*(\Tw^b(R_W)) \cong HH_*(\Tw^b(C_M)).\]
If $W$ has isolated singularities, by Dyckerhoff's generating result and Proposition~\ref{prop:equiv} it follows that $[\Tw^b(C_M)]$ is a compact subcategory of $[\Tw(C_M)]$ (see Section~\ref{sec:generator}). Thus we have an inclusion of dg categories
\[ \Tw^b(C_M) \hookrightarrow \Tw(C_M)^c.\]
Moreover Theorem~\ref{thm:character} implies that every compact object in $\Tw(C_M)$ is a direct factor of an object in $\Tw^b(C_M)$ as $\Psi(\Omega C_M) \in \Tw^b(C_M)$ compactly generates $[\Tw(C_M)]$ by Proposition~\ref{prop:compact_generator}. This implies the above inclusion of categories is an equivalence up to factors, which yields
\[ HH_*(\Tw^b(C_M)) \cong HH_*( \Tw(C_M)^c) \]
by Keller's result~\cite{Keller}. The right hand sided category $\Tw(C_M)^c$ is homotopic to the category $\Tw(\Omega C_M)^c$ via the coproduct preserving homotopy equivalences $\Phi$ and $\Psi$.  As the Hochschild homology is also homotopy invariant, we conclude that
\[ HH_*(\Tw(C_M)^c) \cong HH_*(\Tw(\Omega C_M)^c).\]
Finally the Hochschild homology of $\Tw(\Omega C_M)^c$ can be calculated by that of the dg algebra $\Omega C_M$ by Proposition~\ref{prop:compact_generator}. Combining all these isomorphisms we have shown that
\[ HH_*(\MF(R,W)) \cong HH_*(\Omega C_M).\]
In the following we relate the latter homology group with the Borel-Moore Hochschild homology of the curved algebra $R_W$.

\paragraph{Hochschild homology of $C_M$.} We begin with the classical case where the curvature $W$ is not presented. First we recall the Hochschild homology of a coalgebra $C$. Let $C$ be a coalgebra with a coaugmentation, form the cobar algebra $\Omega C$. The Hochschild chain complex $C_*(C)$ is by definition given by the complex
\[ (\Omega C\otimes^\tau C\otimes^\tau \Omega C) \underset{\Omega C\otimes \Omega C}{\otimes} \Omega C.\]
Here the superscript $\tau$ on tensor symbol is again to denote the twisted tensor product using the natural twisting cochain $\tau: C \rightarrow \Omega C$. Observe that $C_*(C)$ is simply $C\otimes \Omega C$ as a vector space, but the differential is twisted by the natural twisting cochain from $C$ to $\Omega C$. To simply the notations we use $C\tilde{\otimes}\Omega C$ to denote the Hochschild complex $C_*(C)$.

The advantage of this definition of the Hochschild complex for coalgebras is that it is quite simple to relate it to the Hochschild complex of its Koszul dual algebra $\Omega C$. Indeed the latter complex is by definition given by
\[ (\Omega C)\otimes^\tau B\Omega C \otimes^\tau \Omega C)\underset{\Omega C\otimes \Omega C}{\otimes} \Omega C.\]
Notice that these two complexes only differ by the middle term where twisted tensor products are formed. The fact that they are quasi-isomorphic follows from the following classical lemma, see~\cite{Loday} for example.

\begin{lemma}
\label{lem:compo}
Let $C_1\stackrel{\tau_1}{\rightarrow} A$ be a twisting cochain between a dg coalgebra $C_1$ and an dg algebra $A$. Let $C_2\stackrel{\gamma}{\rightarrow} C_1$ be a quasi-isomorphism of dg coalgebras. Then the composition $\tau_2$
\[ C_2 \rightarrow C_1 \rightarrow A\]
is also a twisting cochain. Moreover for any dg $A$-module $F$, the map defined by
\[ C_2\otimes^{\tau_2}F \stackrel{\gamma\otimes \id}{\rightarrow} C_1\otimes^{\tau_1}F\]
is a quasi-isomorphism.
\end{lemma}
We apply the lemma to the unit morphism of the adjunction $\Omega \dashv B$
\[ \eta_C: C\rightarrow B\Omega C\] 
and the natural twisting cochain $C\ra \B\Omega C \ra \Omega C$. The fact the $\eta_C$ is a quasi-isomorphism is well-know for ordinary (dg) algebras (even non-curved $A_\infty$ algebras). We end up with the following quasi-isomorphism between the two Hochschild complexes
\[ C_*(C):=C\tilde{\otimes}\Omega C \stackrel{\eta_C\otimes \id}{\rightarrow} C_*(\Omega C):=B\Omega C \tilde{\otimes}\Omega C.\]
We can add the curvature term $W$ (or $M$) into the previous discussion. All the constructions explained above remain the same as we have already explained the twisting cochain and the twisted tensor products in the curved case in Section~\ref{sec:koszul}. However, the proof of Lemma~\ref{lem:compo} does not generalize as the coalgebra $B\Omega C_M$ is curved with noncommutative coproduct. Hence the differential does not square to zero in this case. It is even problematic to talk about the notion of quasi-isomorphism for these coalgebras. Nevertheless the map $\eta_C\otimes \id$ remains a quasi-isomorphism on the associated Hochschild complexes. This is proved in the following proposition.

\begin{proposition}
\label{prop:hoch}
The map $\eta_C\otimes \id$ is a quasi-isomorphism between the chain complexes $C_*(C_M)$ and $C_*(\Omega C_M)$.
\end{proposition}

\proof Observe the existence of a $\Z$-grading on the space $C_*(C_M)$ by the number of $C$ tensor components. And define the following $\Z$-grading on the space $B\Omega C_M\otimes \Omega C_M$ by
\begin{align*}
\deg(f_1\otimes\cdots\otimes f_k)&:= k \mbox{ for an element in } \Omega C_M;\\
\deg([\alpha_1|\cdots|\alpha_n]\otimes \beta)&:= \deg(\alpha_1)+\cdots+\deg(\alpha_n)+\deg(\beta)-n.
\end{align*}
Then one breaks the Hochschild differentials into two parts. The first part is simply the differential when the curvature is not presented. The second part is the differential defined by the curvature term $M$. For simplicity, we denote them by $d^+$ and $d^-$ respectively. (We will not bother to distinguish them on the two complexes as we will specify the complex when making statements.) Observe that the first differential increases the degrees defined above by $1$ and the second differential decreases the degree by $1$. Hence we have a morphism of mixed complexes
\[ \eta_C\otimes\id :(C_M\tilde{\otimes}\Omega C_M,d^+,d^-) \rightarrow (B\Omega C_M\tilde{\otimes}\Omega C_M,d^+,d^-).\]
Through the associated bi-complex of these mixed complexes (details of the mixed complex technique is explained in~\cite{CT}), we can conclude that the $\eta_C\otimes \id$ is a quasi-isomorphism as it is so on the $E^1$-page. The proof is complete.

\medskip
\remark In the proof it is important that we are dealing with direct sum complexes and $d^+$ is degree increasing, because only in this case the spectral sequences under consideration starts with the differential $d^+$.

\paragraph{Relating to Borel-Moore Hochschild complex.} To relate to the Borel-Moore Hochschild homology we dualize the Hochschild complex $C_*(C_M)$. There is a natural chain map from the Borel-Moore Hochschild chain complex $C_*^{\BM}(R_W)$ of $R_W$ to $C_*(C_M)^\chk$ defined by
\[ R\otimes R^+\cdots R^+\otimes R^+ \hookrightarrow (C\otimes C^+\cdots C^+\otimes C^+)^\chk.\]
This map is in fact a map between mixed complexes whose associated double complexes are isomorphic on the $E^1$-page. This fact follows from classical Hochschild-Konstant-Rosenberg theorem. Strictly speaking the HKR theorem applies only to the left hand side, i.e. for the algebra $R$. But for the right hand side, the Hochschild complex of the coalgebra $C$, it suffice to observe that the Hochschild chain complex $C_*(C)$ is actually double graded by the tensor degree and the polynomial degree. Moreover its graded $k$-linear dual agrees with the Hochschild chain complex of the symmetric algebra $\sym(V^\chk)$ to which we can apply HKR theorem. We summarize the main results obtained in the following theorem.

\begin{theorem}
We have the following isomorphisms:
\begin{align*} 
HH_*(\MF(R,W))\cong &HH_*(\Tw^b(C_M))\cong HH_*(C_M);\\
HH_*^\BM (R_W) &\cong HH_*(C_M)^\chk.
\end{align*}
\end{theorem}

\remark
When $W$ has isolated singularities, the vector space $HH_*(C_M)$ is finite dimensional. Moreover on $HH_*(\Tw^b(R_W))$ there exists a natural non-degenerate pairing  that identifies it with its dual space.

\section{Equivariant matrix factorizations}
\label{sec:orbifold}
In this section we study the orbifold version of Theorem~\ref{thm:koszul} and its applications to categories of equivariant matrix factorizations. Throughout the section we work over the ground field $k=\C$ as we need to consider characters of groups.

\paragraph{Equivariant Koszul duality.} Let $C:=\bS(V)$ to be the symmetric coalgebra over a vector space $V$ and let $M:C\rightarrow k$ be a linear map on $C$ that vanishes on scalar and linear terms. Consider a finite abelian group $G$ acting on $C$ via coalgebra morphisms and that the action preserves the linear map $M$, i. e. the composition
\[ C\stackrel{g}{\rightarrow} C \stackrel{M}{\rightarrow}k \]
is equal to $M$ for any element $g\in G$. Given such data we would like to consider the dg category of equivariant twisted complexes over the curved coalgebra $C_M$. The objects are pairs $(E,Q)$ where $E$ is a cofree $C$-comodule with a $G$-action of the form
\[ E:= \oplus_i C\otimes \C_{\chi_i}.\]
Here $\C_{\chi_i}$ denotes the one dimensional $G$-representation associated to a given character $\chi_i$ and we allow indices to repeat in the direct sum above. The linear map $Q$ is a $G$-equivariant $C$-comodule morphism on $E$. Moreover $Q$ satisfies the matrix cofactorization identity. The morphism spaces between objects would be $G$-equivariant $C$-comodule maps. We denote this category by $\Tw([C_M/G])$ to mimic the orbifold notation. As before we denote by $\Tw^b([C_M/G])$ the full subcategory consisting of finite rank objects. Since the cobar construction is functorial, we also have a $G$-action on the cobar algebra $\Omega C_M$. Thus the category $\Tw([\Omega C_M/G])$ can be defined in a similar way.

The Koszul duality functors $\Phi$ and $\Psi$ are defined in the same way as before. Namely for an equivariant matrix cofactorization $(E,Q)$ define
\[ \Phi(E):= \Omega C_M\otimes ^\tau E\]
where $\Phi(E)$ inherits the tensor product $G$-representation. One can check that the functors $\Phi$ and $\Psi$ send equivariant objects to equivariant objects and equivariant morphisms to equivariant morphisms. Moreover the homotopies constructed in the proof of Theorem~\ref{thm:koszul} can be made $G$-equivariant by averaging if necessary. Thus we arrived the following theorem.

\begin{theorem}
\label{thm:orbifold}
The functors $\Phi$ and $\Psi$ restricted to the equivariant categories to give a homotopy equivalence 
\[\Tw([C_M/G])\cong\Tw([\Omega C_M/G]).\]
\end{theorem}

\paragraph{Smash product algebras.} To make better use of the above Theorem~\ref{thm:orbifold}, we first need to make a change of category. Namely we will switch from equivariant categories to categories of twisted complexes over a smash product algebra. More precisely since $G$ acts on the curved coalgebra $C_M$ in a way that preserves the curved coalgebra structure, we could form the smash product curved coalgebra $C_M\sharp G$. As a vector space it is $C\otimes k[G]$ and the coproduct is defined by
\[ x\otimes g \mapsto \sum_{g_1g_2=g} (x^{(1)} \otimes g_1)\otimes (g_1^{-1}(x^{(2)})\otimes g_2)\]
The curvature of $C_M\sharp G$ is defined by $M$ on the component $C\otimes \id_G$ and zero otherwise. The dg category $\Tw(C_M\sharp G)$ is closely related to the equivariant dg category $\Tw([C_M/ G])$. Observe that the smash product coalgebra $C_M\sharp G$ carry natural $G$-action and $C_M\sharp G$-linear maps are equivalent to $C$-linear maps that are also $G$-equivariant. Thus the category $\Tw(C_M\sharp G)$ is a fully faithful subcategory of $\Tw([C_M/G])$ consists of objects that are free $\C_M\sharp G$-comodules. Conversely every objects of $\Tw([C_M/G])$ is a direct summand of an object in $\Tw(C_M\sharp G)$ through the fully faithful embedding. To see this observe that for any object $(E,Q)\in \Tw([C_M/G])$ form the object
\[ g^*(E,Q):=(\oplus_{g\in G} g^*E, \oplus_{g\in G} g^*Q).\]
One easily checks that $g^*(E,Q)$ is an object of $\Tw(C_M\sharp G)$. Such a relation between the two categories are called equivalence up to factors (from~\cite{Keller}). If two categories are equivalent up to factors, then lots of properties of them are the same. For example (classical) generators of the smaller category are also (classical) generators of the bigger one. It is also proved by Keller~\cite{Keller} that the Hochschild type invariants are isomorphic for these two categories. Observe that $\Phi$ and $\Psi$ restrict to a homotopy equivalence
\[ \Tw(\Omega C_M\sharp G)      \cong \Tw(C_M\sharp G) .\]
As a conclusion we summarize the previous discussion in the following commutative diagram.
\[\begin{CD}
\Tw(\Omega C_M\sharp G) @>\mbox{Koszul duality}>> \Tw(C_M\sharp G) \\
@VV\mbox{inclusion}V                                                      @VV\mbox{inclusion}V \\
\Tw([\Omega C_M/G])    @>\mbox{Koszul duality}>>  \Tw([C_M/G]).
\end{CD}\]
The vertical inclusions are all equivalences up to factors.

\paragraph{Applications to $\MF_G(R,W)$.} The advantage of the smash product construction is that it is clear in this description the object $\Omega C_M\sharp G$ compactly generates the homotopy category of $\Tw(\Omega C_M\sharp G)$. Indeed for an object $F\in \Tw(\Omega C_M\sharp G)$ we have
\[ \Hom_{\Tw(\Omega C_M\sharp G)}(\Omega C_M\sharp G, F)= \Hom_{\Tw(\Omega (C_M))}(\Omega C_M, F)\]
through the inclusion mentioned above. By Corollary~\ref{cor:generator} if the latter is acyclic, then the dg-module $F$ is contractible over $\Omega C_M$. Averaging the contracting homotopy yields a contraction over $\Omega C_M\sharp G$. Hence arguing as in Corollary~\ref{cor:generator} shows that the object $\Omega C_M\sharp G$ compactly generates $[\Tw(\Omega C_M\sharp G)]$. As the categories $\Tw(\Omega C_M\sharp G)$ and $\Tw([\Omega C_M/G])$ are equivalent up to factors, the object $\Omega C_M\sharp G$ (through the inclusion functor) also compactly generates the homotopy category of the latter one.

Applying the Koszul duality functor $\Psi$ yields compact generators for the homotopy category of $\Tw([C_M/G])$. Moreover one can easily identify the generators by observing that the object
$\Omega C_M\sharp G$ when considered as objects in $\Tw([\Omega C_M/G])$ is isomorphic to the direct sum
\[ \oplus_{\chi} \Omega C_M\otimes \C_\chi \]
over the characters of $G$. Hence its image under $\Psi$ is the direct sum
\[ \oplus_\chi \Psi(\Omega C_M)\otimes \C_\chi.\]
Observe that twisting by characters does not change the homology of $\Omega C_M\otimes \C_\chi$ and hence Propositions~\ref{prop:compact_generator}, ~\ref{prop:same} still apply which assert that their $k$-linear duals are (homotopic to) matrix factorizations of the form
\[\left \{ k^\stab \otimes \C_\chi \mbox{ $|$ $\chi$ is a character for the group } G \right \}.\]

\begin{theorem}
\label{thm:orbi_generator}
Let notations be as above and assume that $W$ has isolated singularities. Then the category $[\MF_G(R,W)]$ is classically generated by objects $k^\stab\otimes \C_\chi$.
\end{theorem}

\proof It is enough to show that the subcategory $[\Tw^b(C_M\sharp G)]$ is compact in $[\Tw(C_M\sharp G)]$ in view of Proposition~\ref{prop:equiv}. For this observe that taking cohomology commutes with taking $G$-invariants and hence for a finite rank object $E$ we have
\begin{align*}
\Hom_{[\Tw (C_M\sharp G)]}(E, \oplus E_i) &:= H^0(\Hom_{\Tw (C_M\sharp G)}(E, \oplus E_i)) \\
                                                 &=H^0 (\Hom_{\Tw(C_M)}(E, \oplus E_i))^G\\
                                                 &=[\oplus H^0(\Hom_{\Tw(C_M)} (E, E_i))]^G \\
                                                 &=\oplus \Hom_{[\Tw(C_M\sharp G)]}(E, E_i). 
\end{align*}
Here we have used the fact that $E$ is of finite rank and the group $G$ is finite, which implies that $E$ viewed as an object in $[\Tw(C_M)]$ is compact by Proposition~\ref{prop:equiv}. The theorem is proved.

\paragraph{Equivariant Hochschild homology.} The computation of Hochschild homology of $\MF_G(R,W)$ can be done in the same way as in Section~\ref{sec:hoch}. Again we assume $W$ has isolated singularities throughout the discussion. We begin with an isomorphism
\[ HH_*(\MF_G(R,W)) \cong HH_*(\Tw^b([C_M/G]))\]
as the two dg categories are opposite to each other by the $k$-linear dual functor $D$. Since the compact generators $\Psi(\Omega C_M)\otimes \C_\chi$ of $\Tw([C_M/G])$ lies inside $\Tw^b([C_M/G])$ which is compact under the assumption of $W$ having isolated singularities, we have
\[ HH_*(\Tw^b([C_M/G])) \cong HH_*(\Tw([C_M/G])^c) \cong HH_*(\Tw(C_M\sharp G)^c).\]
The latter isomorphism follows from the fact that the two categories are equivalence up to factors. Finally we invoke the Koszul duality of the curved coalgebra $C_M\sharp G$ which gives a homotopy equivalence
\[  \Tw(C_M\sharp G)^c \cong \Tw(\Omega(C_M\sharp G))^c \]
between dg categories. From this homotopy equivalence and the fact that $\Omega C_M\sharp G$ is a compact generator, we conclude that
\[ HH_*(\Tw(C_M\sharp G)^c) \cong HH_*(\Tw(\Omega(C_M\sharp G))^c) \cong HH_*(\Omega(C_M\sharp G)).\]
Combining the above isomorphisms yields the following isomorphism
\[ HH_*(\MF_G(R,W)) \cong HH_*(\Omega(C_M\sharp G)).\]
Then the same proof as in Section~\ref{sec:hoch} implies the following proposition.
\begin{proposition}
\label{prop:orbihoch}
Let the notations be as above and assume that $W$ has isolated singularities. Then we have the following isomorphisms:
\begin{align*} 
HH_*(\MF_G(R,W)) \cong &HH_*(\Tw^b([C_M/G])) \cong HH_*(C_M\sharp G);\\
HH_*^\BM (R_W\sharp G) &\cong HH_*(C_M\sharp G)^\chk.
\end{align*}
\end{proposition}

\remark The homology groups $HH_*^\BM(R_W\sharp G)$ are explicitly computed in~\cite{CT} via certain localization formula for Borel-Moore homology groups.

\section{Graded matrix factorizations}
\label{sec:graded}
In this section, we study the category of graded matrix factorizations via Koszul duality. The main ideas remain the same as in the orbifold case. The results obtained are closely related to the work of Orlov~\cite{Orlov} (on the relationship between graded matrix factorizations and derived category of coherent sheaves) and Seidel~\cite{Seid} (on the $A_\infty$ category of coherent sheaves on Calabi-Yau hypersurfaces). Throughout the section we work over the ground field $k=\C$.

\paragraph{Gradings.} For a graded commutative ring $S$ and a homogeneous curvature element $W\in R$ of degree $d$, one can define the dg category of graded matrix factorizations $\MF^\gr(R,W)$ (see~\cite{CT} for a definition). As is explained in \textsl{loc. cit.} this category is closely related to certain orbifold construction. We recall some relevant results below.

The symmetric algebra $S:=\bS(V^\chk)$ (non-complete) has a $\Z$-grading by the ordinary polynomial degrees.  The polynomial degree of a homogeneous element $f\in S$ will be denoted by $|f|$. Consider $G:=\Z/d\Z$ acting on $S$ by
\[ \hat{i} (f):= \zeta^{i|f|} f \]
for $\zeta:=\exp(2\pi\sqrt{-1}/d)$, a $d$-th root of unity. Clearly the $G$-action on $S$ preserves the curvature element $W$. This implies that the $G$-action in fact acts on the curved algebra $S_W$. We can then form the smash product curved algebra $S_W\sharp G$. One theorem proved in~\cite{CT} was the fact that graded matrix factorizations can be regarded as $\Z$-graded twisted complexes over $S_W\sharp G$. A subtle point there was that $S_W\sharp G$ does not form a $\Z$-graded curved algebra with the obvious polynomial grading.

To fix this problem we need to introduce a new $\Z$-grading on $S_W\sharp G$. Note that the underlying vector space of $S_W\sharp G$ is $S\otimes k[G]$. The group algebra $k[G]$ has a special basis indexed by characters of $G$. Explicitly we denote by $\chi_i$ for $i\in [0,d-1]$, the characters of the group $G$. They act on $G$ by
\[ \chi_i(\hat{j}):= (\zeta_d)^{i\cdot j} .\]
Then the elements
\[ U_\chi:= \frac{1}{|G|} \sum_{g\in G} \chi(g)\sharp g\]
indexed by these characters form an orthogonal idempotent basis for the group algebra $k[G]$. Using this basis we can define a new $\Z$-grading on the vector space $S\otimes k[G]$. The homogeneous elements are of the form 
\[ f\otimes U_{\chi_j}\]
for some homogeneous polynomial $f\in S$. Define an integer $i\in [0,d-1]$ by
\[ i \equiv j-|f| \pmod d.\]
Then the new grading of $f\otimes U_{\chi_j}$ is defined by
\[ \deg(f\otimes U_{\chi_j}):=  \frac{2}{d}(|f|-j+i).\]
We mention some important properties for this new $\Z$-grading on $S_W\sharp G$. First of all as promised the curvature term $W\sharp \id_G$ has degree $2$ with respect to this grading. To see this observe that
\[ W\sharp \id_G = \sum_{\chi_j} W\otimes U_{\chi_j}.\]
Since $|W|=d$ we have $i=j$ and hence 
\[\deg(W\otimes U_{\chi_j})= \frac{2}{d}\cdot |W|= \frac{2}{d}\cdot d =2 .\] 
Secondly the category of $\Z$-graded twisted complexes over $S_W\sharp G$ is closely related to the category of graded matrix factorizations. In fact it was shown in~\cite{CT} that they are equivalent up to factors. (There we considered $S_W\sharp G$ as a category, then the twist construction would yields in fact an equivalence. Here we prefer to consider $S_W\sharp G$ as a curved algebra.) Namely there is an inclusion 
\[ \Tw_\Z^b(S_W\sharp G) \hookrightarrow \MF^\gr(S,W)\]
which is fully faithful and an equivalence up to factors.

\paragraph{Graded dualization.} Next we dualize the $\Z$-graded curved algebra to consider a $\Z$-graded curved coalgebra $C_M\sharp G$ where $C$ is the symmetric coalgebra $\bS(V)$. We still denote the polynomial degree for a homogeneous $f\in\bS(V)$ by $|f|$. A new $\Z$-grading on $C_M\sharp G$ is defined similarly. Namely homogeneous elements in $C_M\sharp G$ are of the form
\[ f\otimes U_{\chi_j}\]
and the degree of it is given by
\[ \deg(f\otimes U_{\chi_j}):= - \frac{2}{d}(|f|-j+i)\]
for the same $i$ as in the case of algebras. With respect to this $\Z$-grading the map $M:C\rightarrow k$ has degree $2$. Hence it forms a $\Z$-graded curved coalgebra. When forming the category $\Tw^b_\Z(C_M\sharp G)$ we do not want to allow arbitrary arbitrary coalgebra maps but only the direct sums of the homogeneous ones. We introduce a notation to deal with such situations. Let $E$ be a vector space with a $\C^*$-action, we denote by $E^\gr$ the vector space defined by
\[ E^\gr:= \oplus E_i\]
where $E_i$ is the subspace of $E$ on which $\C^*$ acts by $\lambda^i$. With this notation we have
\[ \Hom_{\Tw_\Z(C_M\sharp G)}(-,-):= [\Hom_{\Tw(C_M\sharp G)} (-,-)]^\gr .\]
Then the $\Z$-graded $k$-linear dual operation $D$ defines an equivalence
\[ \Tw^b_\Z(C_M\sharp G)^\oppo \cong \Tw^b_\Z(S_W\sharp G)\]
between dg categories.

\paragraph{$\Z$-graded curved Koszul duality.} 
The next step is to understand the curved Koszul duality for the $\Z$-graded curved coalgebra $C_M\sharp G$. This is easily accomplished by matching the degrees. For this we define a $\Z$-grading on $\Omega C_M\sharp G$ that matches with the new $\Z$-grading on $C_M\sharp G$. The homogeneous elements in $\Omega C_M\sharp G$ are of the form 
\[ [f_1|\cdots|f_k]\otimes U_{\chi_j}\]
for some character $\chi_j$ of the group $G$. Its degree is defined by
\[ \deg([f_1|\cdots|f_k]\otimes U_{\chi_j}):= -\frac{2}{d} (\sum_l |f_l| -j+i)+k \]
where the integer $i\in [0,d-1]$ is defined by 
\[ i \equiv j- \sum_l |f_l| \pmod d .\]
Define the $\Z$-graded Koszul duality functors by (the same formula as before)
\begin{align*}
E \in Tw_\Z(C_M\sharp G) &\stackrel{\Phi}{\mapsto} \Omega(C)\otimes^\tau E \mbox{ and } \\
F \in Tw_\Z(\Omega C_M\sharp G) &\stackrel{\Psi}{\mapsto} C\otimes^\tau F.
\end{align*}
The degrees on $\Phi(E)$ can be defined by
\[ \deg ([f_1|\cdots|f_k]\otimes f_0 \otimes U_{\chi_j}):= -\frac{2}{d} (\sum_{l=0}^k |f_l| -j+i)+k \]
where the integer $i\in [0,d-1]$ is defined by 
\[ i \equiv j- \sum_{l=0}^k |f_l| \pmod d .\]
Similar one can also define degrees for $\Psi(F)$. With respect to these gradings the twisted differentials on $\Phi(E)$ or $\Psi(F)$ have degree one. Moreover it is easy to see that $\Phi$ and $\Psi$ are homotopy equivalences by observing that the homotopy equivalences used in the proof of Theorem~\ref{thm:koszul} respect the new $\Z$-grading (the homotopies are of degree $-1$). 
\begin{theorem}
\label{thm:Z-koszul}
The functors $\Phi$ and $\Psi$ are homotopy inverses between dg categories
\[ \Tw_\Z(\Omega C_M\sharp G) \cong \Tw_\Z(C_M\sharp G).\]
\end{theorem}

\paragraph{Applications to $\MF^\gr(S,W)$.} We assume that $W$ has isolated singularities from now on. One can argue in the same way as in the orbifold case that $\Omega C_M\sharp G$ compactly generates $[\Tw_\Z(\Omega C_M\sharp G)]$. Through the $\Z$-graded Koszul duality functor, $\Psi(\Omega C_M\sharp G)$ defines a compact generator for $[\Tw_\Z(C_M\sharp G)]$. The same proof as in Section~\ref{sec:generator} shows that the object $\Psi(\Omega C_M\sharp G)$ in fact is homotopic to an object in $\Tw^b_\Z(C_M\sharp G)$. Thus its $k$-linear graded dual object in $\MF^\gr(S,W)$ makes sense. To identify this object we consider the natural forgetful functor from $\Tw_\Z(C_M\sharp G)$ to $\Tw (C_M\sharp G)$. Note that this is well-defined as the new $\Z$-grading on $C_M\sharp G$ is in $2\Z$ and hence its reduction modulo $2$ reduces to the purely even grading on the curved coalgebra $C_M\sharp G$. Using the forgetful functor we see that as matrix factorizations the object $D\Psi(\Omega C_M\sharp G)$ is given by
\[ \oplus_i k^\stab \otimes \chi_i.\]
Through the correspondence
\[ \Tw^b_\Z (R_W\sharp G) \hookrightarrow \MF^\gr(S,W)\]
defined in~\cite{CT}, twisting by characters $\chi_j$ corresponds to twisting $(j)$ of ordinary graded $S$-modules. Hence if we assume any lifting of the $\Z$-grading on $k^\stab$, we conclude that the object $D\Psi(\Omega C_M)$ in the category $\MF^\gr(S,W)$ given by the direct sum of the objects
\[ k^\stab(d-1), k^\stab(d-2), \cdots, k^\stab.\]

\begin{theorem}
\label{thm:graded_generators}
Assume that $W$ has isolated singularities, then the collection of objects $k^\stab(d-1), k^\stab(d-2), \cdots, k^\stab$  classically generates $[\MF^\gr(R,W)]$.
\end{theorem}

\proof The theorem follows from the fact that the category $\Tw^b_\Z(C_M\sharp G)$ is compact in $\Tw_\Z(C_M\sharp G)$ which follows from the fact that taking cohomology of a differential of $\Z$-degree $1$ (in particular it is homogeneous) commutes with both taking $G$-invariants and the operation $- \mapsto -^\gr$.

\medskip
\remark In the CY situation, i.e. when $\dim(S)=d=\deg(W)$, the category $[\MF^\gr(S,W)]$ is equivalent to the bounded derived category of coherent sheaves $\D(X)$ on $X:=\Proj S/W$. Denote by $i: X\hookrightarrow \PN^{d-1}:=\Proj S$ the natural embedding of $X$ into the projective space.  Then the above collection of generators corresponds to the collection
\[ i^*\omega_{\PN^{d-1}}[d-1], i^*(\wedge^{d-2}\Omega_{\PN^{d-1}})[d-2],\cdots,i^*\Omega_{\PN^{d-1}}[1],\cO_X\]
through a correspondence $\D(X)\cong [\MF^\gr(S,W)]$. This can be proved by observing that the degree shift in $[\MF^\gr(S,W)]$ corresponds to the composition of the homological degree shift functor and the Seidel-Thomas twist functor associated to the spherical object $\cO_X$ on $\D(X)$, see~\cite{BFK}.

\medskip
\remark The homology of the dg algebra $\Omega C_M\sharp G$ is easily seen to be $\wedge^*(V)\sharp G$. This latter notation is slightly misleading because we did not mean the smash product algebra. It is simply the smash product vector space. The presence of the curvature term puts $A_\infty$ structure on $\wedge^*(V)\sharp G$ via homotopy transfer property. However this computation quickly gets complicated. The author has not been able to describe it even in the case of elliptic curves. We mention two closely related results in these directions. In an unpublished notes~\cite{Seid}, Seidel has obtained the above picture for an $A_\infty$ structure on $\wedge^*(V)\sharp G$ via quite different methods. Explicit calculations for $A_\infty$ structures on elliptic curves have been obtained by Polishchuk in~\cite{Poli}, again through other methods. In latter case even the underlying vector space is different.

\paragraph{Hochschild homology of $\MF^\gr(S,W)$.} The Hochschild homology of the dg category $\MF^\gr(S,W)$ can also be related with the Borel-Moore Hochschild homology of a curved algebra. The proof is the same as the orbifold case except that we use graded $k$-linear dualizing functor. We omit the proof here. The precise results are stated in the following proposition.
\begin{proposition}
\label{prop:graded}
Let the notations be as above and assume that $W$ has isolated singularities. Then we have the following isomorphisms:
\begin{align*} 
HH_*(\MF^\gr(S,W))\cong &HH_*(\Tw^b_\Z([C_M/ G])) \cong HH_*(C_M\sharp G);\\
HH_*^\BM (S_W\sharp G) &\cong HH_*(C_M\sharp G)^\chk
\end{align*}
where the $\chk$ denotes the graded dual operation.
\end{proposition}

\medskip
\remark
Again the groups $HH_*^\BM(S_W\sharp G)$ has been computed in~\cite{CT}. What's new here is the existence of a $\Z$-grading on these homology groups. In the Calabi-Yau situation, the dg version of CY/LG correspondence shows that this computation provides an alternative way to compute the Hochschild homology of CY hypersurfaces.

\appendix
\section{Curved homological perturbation lemma}
\label{sec:chpl}
In this appendix we recall the homological perturbation technique as studied in~\cite{CR}. Then we prove that the homological perturbation lemma remains true when curvatures are presented. This is useful to study homotopy between precomplexes.

In this section we work with a $k$-linear abelian category $\cC$. Our primary application concerns with $\cC$ being the category of $B$-comodules for a coalgebra $B$ over $k$.

\paragraph{Deformation retractions.} Let $(L,b)$ and $(M,d)$ be two complexes over $\cC$. A deformation retraction between them consists of the following data. There are morphisms
\[ i:(L,b) \rightarrow (M,d) \mbox{ and } p:(M,d)\rightarrow (L,b) \]
such that 
\[p\circ i=\id_L.\]
Moreover there is a homotopy $H$ between $i\circ p$ and $\id_M$, i.e. we have
\[ i\circ p= \id+dH+Hd. \]
The triple $(i, p, H)$ is then called a deformation retraction between $(L,b)$ and $(M,d)$. If in additional these maps also satisfy
\begin{equation}
\label{eq:3}
Hi=0,\mbox{ } pH=0,\mbox{ and }  H^2=0,
\end{equation}
Then it is called a special homotopy retraction.

\paragraph{Perturbations.} A perturbation of the complex $(M,d)$ is an odd map $\delta:M\rightarrow M$ such that $(d+\delta)^2=0$. Following the terminologies in~\cite{CR}, we call $\delta$ small if $(\id-\delta H)$ is invertible. For a small perturbation $\delta$, define the operator
\[ A:=(\id-\delta H)^{-1} \delta\]
and define the perturbed homotopy retraction operators by
\begin{equation}
\label{eq:4}
b_1:=b+pAi,\mbox{ } i_1:=i+HAi,\mbox{ }
p_1:=p+pAH,\mbox{ } H_1:=H+HAH.
\end{equation}
Homological perturbation lemma states that the data $(i_1, p_1, H_1)$ defines a new special deformation retraction between the perturbed complexes $(L, b_1)$ and $(M,d+\delta)$. This simple lemma plays an important role in the homotopy theory of algebras.

\paragraph{Curved homological perturbation lemma.} Next we prove a curved version of the homological perturbation lemma. Namely we assume the same initial conditions for $i$, $p$, $H$. But for the perturbation, we do not assume that $(d+\delta)^2=0$. Instead we assume that that
\[ (d+\delta)^2 \mbox{ lies in the center of the algebra } \End(M).\]
We denote this central element by $F:=(d+\delta)^2\in\End(M)$ and call $\delta$ a curved perturbation.

The differential $d_1:=d+\delta$ no longer squares to zero but lies in the center of $\End(M)$.  Such a pair $(M,d_1)$ is called a precomplex. What curved homological perturbation achieves is the fact one can still obtain a deformation retract between precomplexes by perturbing ordinary complexes. The main result of this appendix is the following lemma.

\begin{lemma}
\label{lem:chpl}
(Curved homological perturbation lemma.) Let $(i,p,H)$ a special homotopy retraction data between complexes $(L,b)$ and $(M,d)$. Let $\delta$ be a curved perturbation of $(M,d)$. Then formula~\ref{eq:4} defines a new special homotopy retract between the precomplexes $(L,b_1)$ and $(M,d_1)$ in the following sense:
\begin{itemize}
\item[(A)] $(L,b_1)$ is a precomplex;
\item[(B)] $d_1\circ i_1= i_1\circ b_1$ ($i_1$ is a map of precomplexes);
\item[(C)] $b_1\circ p_1= p_1\circ d_1$ ($p_1$ is a map of precomplexes);
\item[(D)] $p_1\circ i_1=\id_L$ and $i_1\circ p_1=\id_M+d_1H_1+H_1d_1$ (homotopy retract);
\item[(E)] $H_1\circ i_1=0, p_1\circ H_1=0$ and $H_1^2 =0 $ (specialness).
\end{itemize}
\end{lemma}

\proof The proof is analogous to the proof of the ordinary perturbation lemma in~\cite{CR}. We basically only need to check the above formulas with a weaker condition that $F$ is in the center (weaker as $0$ is in the center). We begin with the following lemma.

\begin{lemma}
We have
\begin{itemize}
\item $\delta HA=AH\delta =A-\delta$;
\item $(\id-\delta H)^{-1}=\id+AH$ and $(\id-H\delta)^{-1}=\id+HA$;
\item $AipA+Ad+dA=F+FAH+FHA$.
\end{itemize}
\end{lemma}

\proof The first two equations are direct computations and is the same as in~\cite{CR}. For the last one, we have
\begin{align*}
 AipA+Ad+dA &= A(\id+dH+Hd)A +Ad +dA \\
            &= A^2 + AdHA +AHdA +Ad +dA \\
            &= A^2 + Ad(HA+\id)+(AH+\id)dA \\
            &= A^2 + Ad(\id-H\delta)^{-1}+(\id-\delta H)^{-1}dA\\
            &= (\id-\delta H)^{-1}[(\id-\delta H)A^2(\id-H\delta)+\\
&    +(\id-\delta H)Ad+dA(\id-H\delta)](\id-H\delta)^{-1}\\
            &= (\id-\delta H)^{-1}[\delta^2+\delta d +d \delta] (\id-H\delta){-1}\\
            &= F(\id-\delta H)^{-1}(\id-H\delta)^{-1}\\
            &= F(\id+AH)(\id+HA)\\
            &= F+FAH+FHA.
\end{align*}
With these preparations, the proof of Lemma~\ref{lem:chpl} follows easily as an extension of the case without curvature. Let us first prove part $(A)$.
\begin{align*}
b_1^2 &= (b+pAi)(b+pAi) \\
      &= bpAi+pAib+p(AipA)i \\
      &= bpAi+pAib+p(F+FAH+FHA-Ad-dA)i \\
      &= pFi+pFAHi+pFHAi\\
      &= pFi+pFAHi+pHFAi \mbox{ ($F$ is central) }\\
      &= pFi \mbox{ (specialness) }.
\end{align*}
Thus $b_1^2$ is simply the restriction of $F$ on its subspace $L$ (via $i$ and $p$). And hence it is in the center of $\End(L)$, which proves that $(L,b_1)$ is a precomplex. For part $(B)$ we have
\begin{align*}
i_1b_1-(d+\delta)i_1 &=(i+HAi)(b+pAi)-(d+\delta)(i+HAi)\\
                     &=ib+ipAi+HAib+H(AipA)i-di-dHAi-\delta i- \delta HAi \\
                     &=ipAi+HAib+H(F+FAH+FHA-dA-Ad)i -\\
&-dHAi-\delta i-(A-\delta)i\\
                     &=ipAi-HdAi-dHAi-Ai+HFi+HFAHi+HFHAi\\
                     &=(ip-Hd-dH-\id)Ai+HFi+HFAHi+HFHAi\\
                     &=FHi+HFAHi+FHHAi\\
                     &=0\mbox{ (by specialness and $F$ is central) }.
\end{align*}
Similarly we check that $p_1$ is map of precomplexes:
\begin{align*}
b_1p_1-p_1(d+\delta)&= (b+pAi)(p+pAH)-(p+pAH)(d+\delta)\\
                    &= bpAH+pAip+p(AipA)H-p\delta-pAHd-p(AH\delta)\\
                    &= bpAH+pAip-p(Ad+dA)H+p(F+FAH+FHA)H \\
                    &-p\delta-pAHd-p(A-\delta)\\
                    &= pAip-pAdH-pAHd-pA+pFH+pFAHH+pFHAH                    \\
                    &= pFH+pFAHH+pFHAH\\
                    &=0 \mbox{ (by specialness and $F$ is central) }.
\end{align*}
This proves part $(C)$. For part $(D)$ we have
\begin{align*}
p_1i_1 &=(p+pAH)(i+HAi)\\
       &=pi+pHAi+pAHi+pAHHAi\\
       &=\id \mbox{ (by specialness) }.
\end{align*}
In the reversed direction, we have
\begin{align*}
\id+H_1d_1& +d_1H_1-i_1p_1= \id+(H+HAH)(d+\delta)+\\
&+(d+\delta)(H+HAH)-(i+HAi)(p+pAH)\\
                         &= H\delta+HAHd+H(AH\delta)+\delta H+dHAd+ (\delta HA)H \\
                         &- ipAH-HAip-H(AipA)H\\
                         &= H\delta+HAHd+H(A-\delta)+\delta H+dHAd+ (A-\delta)H \\
                         &- ipAH-HAip+H(Ad+dA)H-H(F+FAH+FHA)H\\
                         &= HA(Hd+\id+dH-ip)+(dH+\id-ip+Hd)AH\\
                         &-HFH-HFAHH-HFHAH\\
                         &=0 \mbox{ (again by specialness and $F$ is central) }.
\end{align*}
Thus we have shown that $(i_1, p_1, H_1)$ forms a deformation retraction. It still remains to prove part $(E)$. This is again a computation.
\begin{align*}
H_1\circ i_1 &= (H+HAH)(i+HAi)\\
             &= Hi+HHAi+HAHi+HAHHAi\\
             &= 0;
\end{align*}
\begin{align*}
p_1\circ H_1 &=(p+pAH)(H+HAH)\\
             &=pH+pHAH+pAHH+pAHHAH\\
             &=0;
\end{align*}
\begin{align*}
H_1\circ H_1 &= (H+HAH)(H+HAH)\\
             &=HH+HAHH+HHAH+HAHHAH\\
             &=0.
\end{align*}
Thus the lemma is proved.

\end{document}